\documentclass[oneside,frenchb,english,reqno]{amsart}
\usepackage{palatino}
\usepackage[T1]{fontenc}
\usepackage[latin1]{inputenc}
\usepackage{amssymb}

\makeatletter

\providecommand{\LyX}{L\kern-.1667em\lower.25em\hbox{Y}\kern-.125emX\@}

\theoremstyle{plain}
\newtheorem*{lem*}{Lemma}
\theoremstyle{plain}
\newtheorem{thm}{Theorem}

\usepackage{babel}
\makeatother
\begin{document}

\title{Maximal smoothness of the anti-analytic part of a trigonometric null
series}

\author{Gady Kozma }
\email{gadykozma@hotmail.com}

\author{Alexander Olevski\u\i{}}
\email{olevskii@post.tau.ac.il}

\address{Tel Aviv University, School of Mathematics, Ramat Aviv 69978, Israel}

\thanks{Research supported in part by the Israel Science Foundation}

\begin{abstract}
We proved recently \cite{KO03} that the anti-analytic part of a trigonometric
series ,converging to zero almost everywhere, may belong to $L^{2}$ on the circle.
 Here we prove that it
can even be $C^{\infty }$, and we characterize precisely the possible
degree of smoothness in terms of the rate of decrease of the Fourier
coefficients. This sharp condition might be viewed as a {}``new 
quasi-analyticity ''.
\end{abstract}
\maketitle

\section{RESULTS}

\markright{Maximal smoothness of anti-analytic null series}

The classical Menshov example shows that a (nontrivial) trigonometric
series\begin{equation}
\sum c(n)e^{int}\label{eq:sumcneint}\end{equation}
may converge to zero almost everywhere (a.e.). Such a series is called
a null series. This result was the origin of modern uniqueness theory
in Fourier Analysis, see \cite{B64,KS94,KL87}. A null series can
not be analytic, that is involve positive frequencies only. This follows
from Abel and Privalov theorems. On the other hand, we proved recently
\cite{KO03} that the anti-analytic part can be small in the sense
that\begin{equation}
\sum _{n<0}|c(n)|^{2}<\infty .\label{eq:sumcn2fin}\end{equation}
It turns out that a much stronger property is possible: the anti-analytic
part can be infinitely smooth.

\begin{thm}
\label{thm:cinf}There exists a trigonomteric series (\ref{eq:sumcneint})
convergent to zero a.e., such that\[
c(n)=O(1/|n|^{k})\; (n<0)\textrm{ for every }k=1,2,\dotsc \]

\end{thm}
Moreover the following result is true :

\begin{thm}
\label{thm:precise}Let $\omega $ be a function $\mathbb{R}^{+}\to 
\mathbb{R}^{+}$,
$\omega (t)/t$ concave and\begin{equation}
\sum \frac{1}{\omega (n)}=\infty .\label{eq:sumwinf}\end{equation}
Then there exists a null-series such that the amplitudes in the negative
spectrum satisfy the condition:\begin{equation}
c(n)=O(\exp (-\omega (\log |n|))),\; n<0.\label{eq:qa}\end{equation}

\end{thm}
It is remarkable that the condition is sharp. The following uniqueness
theorem is true.

\begin{thm}
\label{thm:unique}If a series (\ref{eq:sumcneint}) converges to
zero a.e., and the coefficients satisfy the condition (\ref{eq:qa}),
where $\omega (t)/t$ increase and \begin{equation}
\sum \frac{1}{\omega (n)}<\infty ,\label{eq:sumwfin}\end{equation}
then $c(n)=0$ for all $n\in \mathbb{Z}$.
\end{thm}
So for series (\ref{eq:sumcneint}) converging a.e.~on the circle,
(\ref{eq:qa}) and (\ref{eq:sumwfin}) appears as a sharp quasi-analyticity
condition for the amplitudes of the negative spectrum, which ensures
the uniqueness property. Those amplitudes of a null series may, for
example, decrease as $n^{-\log \log n}$ but not as $n^{-(\log \log n)^2}$.

With respect to Theorem \ref{thm:unique}, it should be mentioned
that if one replaces convergence a.e.~with convergence on a set $E$
of positive mesaure, then a sharp uniqueness condition is the usual
quasi-analyticity:\begin{equation}
c(n)=O(\exp (-\rho (n)))\quad (n<0),\quad \sum \frac{\rho 
(n)}{n^{2}}=\infty .\label{eq:quasianal}\end{equation}
This follows from Beurling theorem \cite{B} extended by Borichev
\cite{B88}, which implies that a series (\ref{eq:sumcneint}), 
(\ref{eq:quasianal})
converging on $E$ to zero is trivial. The sharpness follows from classical
results, see \cite{M35}. In fact, in \cite{B88}, the sum of the analytic part
of (\ref{eq:sumcneint})
is understood (like in Privalov theorem) as a non-tangential boundary limit,
which is assumed to exist on $E$.  In this setting uniqueness holds under doubly
exponentional growth condition of this part in the disc.
 Our theorem 
\ref{thm:unique}
also admits such a version, but the growth conditions necessary are
much stronger.

Below we give a sketch of the ideas involved in the proof of theorem \ref{thm:cinf}.
Theorem \ref{thm:precise} can be obtained basically by the same approach.
We do not discuss here the proof of theorem \ref{thm:unique}.

In the proof below we construct a probabilistically-skewed ``thick'' Cantor
set $K$ of measure zero and a  random harmonic function
$f$ on the disk with singularities on $K$.
 Taking $F=\exp 
(f+i\widetilde{f})$
and denoting by $F^{*}$ the boundary value of $F$ on the circle,
we shall show that $F^{*}$ is smooth and that the Taylor coefficients
$\widehat{F}(n)\to 0$ with probability $1$. Hence the coefficients
$c(n):=\widehat{F}(n)-\widehat{F^{*}}(n)$ are the Fourier coefficients
of a singular compactly supported distribution on $\mathbb{T}$ and
$c(n)\to 0$, which gives, by \cite[p.\ 54]{KS94}, that (\ref{eq:sumcneint})
converges to zero almost everywhere, as required.

It is interesting to compare the proof to the one used in \cite{KO03}.
There $f$ was the Poisson integral of a singular (non-stochastic)
measure on $K$. This approach, however, cannot work here, even if
$f$ is taken to be the sum of a singular measure and an $L^{1}$
function.

\section{\label{sec:Construction}Construction}

Let \begin{align}
\sigma _{n} & :=\frac{1}{2^{n}\log (n+2)}\quad n\geq 1,\quad \sigma 
_{0}=1,\label{eq:defsigman}\\
\tau _{n}  &:=\frac{1}{12}(\sigma _{n-1}-2\sigma _{n})\approx 
\frac{1}{2^{n}n\log ^{2}n}\label{eq:deftaun}
\end{align}
where $X\approx Y$ stands, as usual, for $cX\leq Y\leq CX$, and where $c$
and $C$ stand, here and everywhere, for some absolute constants. Let
$l\in C^{\infty }\left]0,1\right]$ 
be a function satisfying 
$l(x)  =-\log ^{2}x$ for $x  <1/3$ , 
$l(x)  =-1$ for $x > 2/3$ and $l\leq -1$ everywhere.
Given $s \in [0,1]$ define functions on $\mathbb{R}$
\[
l^{\pm }(x;s):=\begin{cases}
l(x) & 0<x\leq 1\\
-1 & 1<x\leq 2\pm s\\
l(3\pm s-x) & 2\pm s<x\leq 3\pm s\end{cases}\]
and $0$ otherwise.

Assume at the $n$'th step of induction that we have $2^{n}$ intervals $I(n,k)$
of length $\sigma _{n}$ (intervals of rank
$n$), and let $K_{n}:=\bigcup _{k=0}^{2^{n}-1}I(n,k)$;  assume
also we have a function $f_{n}:[0,1]\to \mathbb{R}$ such that 
$f_{n}|_{I(n,k)}=M(n)$
i.e.~some constant independent of $k$. Examine one $I=I(n,k)$.
Divide $I$ into two equal parts, $I=I'\cup I''$.
$f_{n+1}$ will now be defined on the sides of $I'$ using
\begin{equation}
f_{n+1}:=\begin{cases}
n\sqrt{\log n}\cdot l^{+}(x/\tau_{n+1};s)&\mbox{left side of $I'$}\\
n\sqrt{\log n}\cdot l^{-}(x/\tau_{n+1};s)&\mbox{right side of $I'$}
\end{cases}
\label{eq:deffnp1}
\end{equation}
which leaves a space of $\frac{1}{2}\sigma_n - 6\tau_{n+1}=\sigma_{n+1}$
in $I'$ undefined --- this will be $I(n+1,2k)$. We fix
$M_{n+1}$ from the condition\begin{equation}
\int _{I'}f_{n+1}-f_{n}=0\label{eq:intIzero}\end{equation}
and it is clear that $M_{n+1}$ does not depend on $s$. Repeat the
construction inside $I''$ with $s=s(n+1,2k+1)$. We remark that the factor
$n\sqrt{\log n}$ in (\ref{eq:deffnp1}), or to be more precise, the fact
that it is superlinear, is the one that guarantees that
the final function $F$ is $C^\infty$.

For now the choice of the $s(n,k)$ is arbitrary. It is only for the last step, 
that we will take the $s(n,k)$ to be random (independent and uniformly
distributed on $[0,1]$). Then we will prove that a
null series with smooth anti-analytic part is generated for almost
any choice of $s$-es.

\section{Estimates}

\subsection{The maximum of $f_n$}

The magnitude of the $\tau_{n}$ (\ref{eq:deftaun}) together with (\ref{eq:deffnp1}) gives
that the negative part of  $f_{n}$  of rank  $i$  has integral
$\approx\log^{-3/2}i$, and hence a sum and (\ref{eq:defsigman}) gives

\begin{equation}
  M_{n}\approx \frac{n}{\sqrt{\log n}}.
  \label{eq:maxfn}
\end{equation}

Similarly, for any interval $I$ of rank $n-1$, \begin{equation}
\int _{I}|f_{n}(x)-f_{n-1}(x)|\stackrel{(\ref {eq:intIzero})}{=}2\int 
_{I}(f_{n}-f_{n-1})^{-}\approx \frac{1}{2^{n}\log ^{3/2}n}\leq 
C2^{-n}.\label{eq:fnfn12n}\end{equation}
We remark that the fact that $M_n$ is sublinear is the one that guarantees
that our final $F$ will have $\widehat{F}(m)\to 0$. Hence the proof hinges
around the following, somewhat paradoxical situation: even though $K$ has
measure zero, it is sufficiently thick so that it would be possible to
balance superlinear growth outside $K$ (the $n\sqrt{\log n}$ factor in
(\ref{eq:deffnp1})) with sublinear growth inside $K$. The proof of theorem
\ref{thm:precise} explores this effect to its maximum.

\subsection{The limit of the $f_n$}

We identify $[0,1]$ with the circle $\{|z|=1\}$ , extend $f_{n}$
as harmonic functions into the disk $\mathbb{D}$ and denote the extensions
by $f_{n}$ as well. We need to estimate $f_{n}$ and their derivatives $f^D_{n}$
( we mean  tangential derivative, i.e.~if $f=f(re^{2\pi 
i\theta })$
then $f':=\frac{df}{d\theta }$).  Using (\ref{eq:intIzero}),
(\ref{eq:fnfn12n}), integration by parts and standard estimates for
the derivatives of the Poisson kernel one can prove:\begin{equation}
|f_{n+1}^{(D)}(z)-f_{n}^{(D)}(z)|\leq 
\frac{C(D)}{2^{n}d(z,K_{n})^{D+1}}\quad \forall z\in 
\overline{\mathbb{D}}\setminus K_{n},\; \forall D\in \{0,1,\dotsc 
\}\label{eq:fnDfn1DdzK}\end{equation}
where $d(z,K)$ denotes the distance of the point $z$ from the set
$K$. Denote by $\widetilde{f_{n}}$ the harmonic conjugate of $f_{n}$.
Using the conjugate Poisson kernel we get the same estimate for 
$|\widetilde{f_{n+1}}^{(D)}(z)-\widetilde{f_{n}}^{(D)}(z)|$.

These two inequalities show that $f_{n}$ and $\widetilde{f_{n}}$ converge
uniformly on compact subsets of $\overline{\mathbb{D}}\setminus K$. Denote
their limits by $f$ and $\widetilde{f}$ respectively --- $\lim \widetilde{f_{n}}$ is
clearly the conjugate of $\lim f_{n}$, which justifies the notation
$\widetilde{f}$.

The boundary values of $f$ are simple to estimate, as $f|_{[0,1]\setminus 
K_{n}}\equiv f_{n}|_{[0,1]\setminus K_{n}}$.
Hence, directly from the definitions of $f_{n}$ and $l$ we get that
$f$ has singularities on $K$ and on a countable set of points $Q$
--- the boundaries and middles of all the intervals $I(n,k)$. Denote
$K':=K\cup Q$. From (\ref{eq:fnDfn1DdzK}) and properties of $f_n$ one can
deduce that on $\mathbb{T}$,
\begin{equation}
|f^{(D)}(x)|\leq \frac{C(D)}{d(x,K')^{D+1}}.\label{eq:fdrough}\end{equation}
This also holds for $\widetilde{f}^{(D)}$, though it
is necessary to first prove an analog of (\ref{eq:fdrough})
for $\widetilde{f_{n}}$ uniformly in $n$ and take the limit as $n\to \infty 
$.
The estimate for $\widetilde{f_{n}}$ follows in turn from the estimate
for $f_{n}$ and  estimates on the derivatives of the Hilbert
kernel.

\subsection{Smoothness}

Define now $F=\exp (f+i\widetilde{f})$. We  use the notation
$F^{*}$ for  the boundary value, considered as a function on
$\mathbb{T}$, in order to distinguish it from the {}``true'' limit
value of $F$ on the boundary of the circle which is a distribution
with a singular part supported on $K$. We note that $F$ is not in
$H^{\infty }$  and therefore the coefficients
$c(n)=\widehat{F}(n)-\widehat{F^{*}}(n)$ are non-trivial.

A rather straightforward calculation starting from (\ref{eq:maxfn})
shows that
\begin{equation}
f(x)\leq -c\log \frac{1}{d(x,K')}\sqrt{\log \log \frac{1}{d(x,K')}}\quad 
\forall x\in \mathbb{T}\setminus K'.
\end{equation}
Combining the fact that $f$ goes to $-\infty$ faster than $\log 1/d(z,K')$
with the rough estimates of (\ref{eq:fdrough}) (and the corresponding
inequality for $\widetilde{f}$) one can prove that $F^{*}\in C^{\infty }([0,1])$.

\section{\label{sec:Fmhatzero} Probability}

Denote $F_{n}=\exp (f_{n}+i\widetilde{f_{n}})$ for
$n=n(m)=\left\lfloor C\log m\right\rfloor $.
Then another relatively simple conclusion from (\ref{eq:fnDfn1DdzK}) is
that for some $C$ sufficiently
large, the following inequality for Taylor coefficients holds\begin{equation}
|\widehat{F_{n}}(m)-\widehat{F}(m)|=\int 
_{(1-1/m)\mathbb{T}}z^{-m-1}(F_{n}(z)-F(z))\, dz\leq 
\frac{C}{m}.\label{eq:FnmFm}\end{equation}

We shall not give many details for the probabilistic argument. In general
it uses a fourth moment calculation. Define therefore, for every
$0\leq k<2^{n}$,
\[
\mathcal{I}_{k}=\int _{I(n,k)}F_{n}(x)e^{imx}\, dx\]
For which we have an absolute bound (from (\ref{eq:maxfn}))
\begin{equation}
|\mathcal{I}_{k}|\leq \int _{I(n,k)}|F_{n}(x)|\leq \sigma 
_{n}e^{Cn/\sqrt{\log n}}=:\gamma_n=\gamma.
\end{equation}

\begin{lem*}
\label{lem:complex}Let $0\leq k_{1},k_{2},k_{3},k_{4}<2^{n}$ and
let $1\leq r\leq n$, and assume that $I(n,k_{i})$ belong to at least
three different intervals of rank $r$. Then \[
\mathbb{E}(\mathcal{I}_{k_{1}}\mathcal{I}_{k_{2}}\mathcal{I}_{k_{3}}\mathcal{I}_{k_{4}})\leq 
\gamma ^{4}\frac{C\log ^{4}m}{m^{2}\tau _{r}^{3}}.\]

\end{lem*}

Had we needed to estimate $\int f_n(x)e^{imx}$ the lemma would have been
standard, since $f$ has a local structure and by conditioning on the location
of the intervals of rank $r$ we would achieve independence between the various
$\mathcal{I}_k$-s. However, $F$ contains also the $\widetilde{f}$ component
which is non-local. Still, it turns out that after the conditioning step we are
left with a function of two variables which can be estimated by two (rather
long) integrations by parts. We skip this calculation entirely.

Proceeding with the proof of the theorem, define \[
X=X_{m}=\sum _{k=0}^{2^{n}-1}\int _{I(n,k)}F_{n}(x)e^{imx}\, dx.\]
The difference $\widehat{F_{n}}(m)-X$ is the integral over the subset of
$\mathbb{T}$ where $f_n=f$, and there $F_n$ is $C^{\infty }$ uniformly in $n$, and in
particular this integral is $\leq C/m$. Therefore we want to bound $X$, and we
shall estimate $\mathbb{E}X^{4}$. Let\[
E(k_{1},k_{2},k_{3},k_{4}):=\mathbb{E}\prod \mathcal{I}_{k_{i}}\]
let $r(k_{1},\dotsc ,k_{4})$ be the minimal $r$ such that $I(n,k_{i})$
are contained in at least $3$ different intervals of rank $r$. A
simple calculation shows \[
\#\{(k_{1},\dotsc ,k_{4}):r(k_{1},\dotsc ,k_{4})=r\}\approx 2^{4n-2r}.\]
The estimate of the lemma is useless if $r$ is too
large. Let $R$ be some number. For $r\geq R$ use the simple $|E(k_{1},\dotsc 
,k_{4})|\leq \gamma ^{4}$
and $\gamma=2^{-n}m^{o(1)}$ to get\begin{equation}
E_{1}:=\sum _{r(k_{1},\dotsc ,k_{4})\geq R}E(k_{1},\dotsc ,k_{4})\leq 
C\gamma ^{4}2^{4n-2R}\leq m^{o(1)}2^{-2R}.\label{eq:E1}\end{equation}
For smaller $r$, we use the lemma to get
$E(k_{1},\dotsc ,k_{4})\leq \gamma ^{4}m^{-2+o(1)}\tau _{r}^{-3}$
and then, using $\tau _{r}=2^{-r+o(r)}$,\begin{align}
E_{2}: & =\sum _{r(k_{1},\dotsc ,k_{4})<R}E(k_{1},\dotsc ,k_{4})\leq \gamma 
^{4}2^{4n}m^{-2+o(1)}\sum _{r=1}^{R}2^{-2r}\tau _{r}^{-3}=\nonumber \\
& =m^{-2+o(1)}\sum _{r=1}^{R}2^{r+o(r)}=m^{-2+o(1)}2^{R+o(R)}.
\end{align}
Picking $R=\left\lfloor \frac{2}{3}\log m\right\rfloor $ we get from
(14) and (15) that $\mathbb{E}X_{m}^{4}\leq m^{-4/3+o(1)}$
and hence $\mathbb{E}\left(\sum X_{m}^{4}\right)<\infty $ and in
particular $X_{m}^{4}\to 0$ with probability 1. As remarked above,
this shows that $\widehat{F_{n}}(m)\rightarrow 0$ and hence using
(\ref{eq:FnmFm}) that $\widehat{F}(m)\to 0$ .\qed

\end{document}